\documentclass[11pt]{article}

\usepackage{graphics,enumitem,epsfig,textcomp}
\usepackage{amsfonts,amsmath,amssymb,amsthm}
\usepackage{euscript,color,mathrsfs}

\usepackage{booktabs}

\usepackage{url}
\usepackage{hyperref}
\hypersetup{
    colorlinks=true,
    linkcolor=blue,
    filecolor=blue,      
    urlcolor=blue,
    citecolor=blue,
    }

\usepackage[mathlines]{lineno}

\usepackage[margin=1in]{geometry}

\usepackage{authblk}

\numberwithin{equation}{section}

\newtheorem{proposition}{Proposition}
\newtheorem{remark}{Remark}

\def\begproof{\noindent{\bf Proof: }}
\def\endproof{\par\rightline{\vrule height5pt width5pt depth0pt}\medskip}
\def\div{\nabla\cdot}

\def\d{\,\mathrm{d}}
\def\N{\mathbb{N}}
\def\Z{\mathbb{Z}}
\def\R{\mathbb{R}}

\def\({\begin{eqnarray}}
\def\){\end{eqnarray}}
\def\[{\begin{eqnarray*}}
\def\]{\end{eqnarray*}}
\def\part#1#2{\frac{\partial #1}{\partial #2}}

\def\grad{\nabla}
\def\Norm#1{\left\| #1 \right\|}

\def\tot#1#2{\frac{\d #1}{\d #2}} 
\def\totk#1#2#3{{\frac{\d^#3 #1}{\d #2^#3}}}
\def\laplace{\Delta}
\def\d{\,\mathrm{d}}
\def\N{\mathbb{N}}
\def\R{\mathbb{R}}

\def\epsilon{\varepsilon}

\def\d{\mathrm{d}}

\def\x{x}

\def\z{z}

\def\mm{\mathfrak{m}}

\def\Td{\mathbb T^d}

\title{Gradient Flow Structure of the Spontaneous Aggregation Model}
\date{}

\author{
     Jan Haskovec\\ Mathematical and Computer Sciences
            and Engineering Division,\\
         King Abdullah University of Science and Technology,\\
         Thuwal 23955-6900, Kingdom of Saudi Arabia\\
         {\it jan.haskovec@kaust.edu.sa}
     }

\begin{document}

\maketitle

\begin{abstract} \noindent
We identify a previously unnoticed gradient-flow structure of the
Fokker--Planck equation arising in the spontaneous particle aggregation model.
For an exponential response function, the equation is a generalized
Wasserstein gradient flow of the attractive McKean--Vlasov free energy
with a nonlocal mobility. We show that, within a natural class of local
entropies and symmetric interaction energies, this structure essentially
singles out the exponential response. We discuss consequences for
stationary states, global minimizers, the single-cluster structure in one
dimension, and formally outline convergence to equilibrium using a
Wasserstein--\mbox{\L}ojasiewicz inequality.
\end{abstract}
\vskip 5mm

\noindent
\textbf{Keywords.} Spontaneous aggregation; Fokker--Planck equation; Gradient flow.
\vskip 5mm


\section{Introduction}\label{sec:Intro}
Aggregation phenomena arise in a wide range of biological systems, where
individuals form spatial clusters through local interactions with their
neighbors. In many classical models this behavior is generated by an
explicit attractive force or a directed drift towards regions of high
density. The spontaneous aggregation model introduced in
\cite{BHW:2012}, originally under the name ``direct aggregation'', is based
on a different mechanism. Individuals move without a preferred direction,
but their diffusivity decreases with the locally perceived population
density. Consequently, particles tend to remain in densely populated
regions, providing a positive feedback mechanism that may lead to the
spontaneous formation of clusters. Its formal mean-field limit is the
nonlinear and nonlocal Fokker--Planck equation
\(   \label{eq:FP}
   \part{\varrho}{t} =\frac12\Delta\bigl(G(W\ast\varrho)^2\varrho\bigr),
\)
where the response function $G$ describes the dependence of the motility on
the locally averaged density $W\ast\varrho$. The model exhibits aggregation,
pattern formation and slow coarsening in suitable parameter regimes; see
\cite{BHW:2012,EH:2026,HO:2015}.

The main purpose of the present note is to point out a structural property
of this equation that, to the best of the author's knowledge, has not been
reported previously. Namely, if the response function is exponential,
\(  \label{eq:Gexp}
   G(s)=\exp(-\beta s),\qquad \beta>0,
\)
then the Fokker--Planck equation can be written as a generalized gradient
flow of the attractive free energy
\[
   \mathcal F_\beta(\varrho) = \int_{\Td}\varrho \, (\log\varrho-1)\,\d x  -\beta\int_{\Td}\varrho(W\ast\varrho)\,\d x,
\]
with the nonlocal mobility
\[
   \mathfrak m(\varrho)
   =\frac12\varrho\exp(-2\beta W\ast\varrho).
\]
This is the same free energy as for the attractive McKean--Vlasov equation
\cite{Carrillo-Gvalani-Pavliotis-Schlichting:2020, Carrillo-McCann-Villani:2003},
but the two models possess different Onsager operators. They therefore
share the same mass constrained critical points, stationary states and global
minimizers, while their transient dynamics and convergence rates may be
substantially different.

It is particularly interesting that the exponential response function \eqref{eq:Gexp}
appears to be singled out by the gradient-flow requirement. Within the
natural class of free energies composed of a convex local entropy and a
symmetric quadratic interaction energy, comparison of the fluxes forces the
local entropy to be, up to affine terms and scaling, the Boltzmann entropy,
and simultaneously forces the logarithmic derivative $G'/G$ to be constant.
Thus $G$ must be exponential. The same conclusion follows from the symmetry
condition for the Fr\'echet derivative of the natural candidate chemical
potential. This does not constitute an impossibility result for every
conceivable variational formulation: we do not exclude the possibility
that a nonexponential response function might admit a gradient structure
based on a different metric or on additional state variables.
Nevertheless, within the broad and natural
class considered here, the exponential response is necessary.

The gradient-flow formulation provides a unified way of deriving several
properties of the model. We characterize positive stationary states as
constrained critical points of the free energy and obtain their Gibbs-type
self-consistency equation and quantitative bounds. We relate the loss of
local minimality of the homogeneous state to its linear instability,
establish existence and elementary uniqueness criteria for global
minimizers, and use the periodic Riesz rearrangement inequality to show that
in one spatial dimension a global minimizer can be chosen to consist of a
single symmetric cluster. Finally, we explain formally how the recent
Wasserstein--\mbox{\L}ojasiewicz inequality for McKean--Vlasov energies
\cite{Choi:2026} may be combined with the nonlocal mobility to obtain
convergence of bounded trajectories to individual stationary states.

The present work is intended as a short structural note rather than a
complete analytical treatment. We therefore emphasize the main ideas and
their consequences without developing a full theory of well-posedness,
regularity and asymptotic compactness. In particular, the discussion of
long-time convergence is formal and indicates the assumptions and arguments
that would be required for a rigorous result. A detailed analysis of these
questions is left for future work.

The paper is organized as follows. In Section~\ref{sec:TheModel} we formulate
the main assumptions of the model
and recall the linearized stability criterion for its homogeneous steady states.
Section~\ref{sec:GradientFlow} derives the gradient-flow structure, explains
why it selects the exponential response function, and discusses its relation
to the classical McKean--Vlasov equation and the variational properties of
the free energy. In Section~\ref{sec:StationaryStates} we characterize
stationary solutions and study the homogeneous state and the global
minimizers. Section~\ref{sec:SingleCluster} establishes the single-cluster
structure of global minimizers in one dimension. Finally,
Section~\ref{sec:LongTime} gives a formal outline of the long-time convergence
argument based on the Wasserstein--\mbox{\L}ojasiewicz inequality.

\section{The model and main assumptions}\label{sec:TheModel}

The mean-field of the spontaneous particle aggregation model derived  in~\cite{BHW:2012} 
consists of the Fokker-Planck equation \eqref{eq:FP} for the particle number density $\varrho = \varrho(t,{\x})$.
We pose it on the $d$-dimensional torus $\Td$.
Throughout the paper we assume the sensing kernel $W=W({\x})$ to be nonnegative, bounded,
radially symmetric and nonincreasing with respect to $|\x|$.
We regard $W$ as a compactly supported function on $\mathbb R^d$, periodically extended to the torus $\Td$, and normalized by
\begin{equation}  
\label{eq:Wnorm}
\int_{\Td} W({\x}) \, \d {\x} = 1 \,. 
\end{equation}
The convolution in \eqref{eq:FP} is given by
\[
   W\ast \varrho(t, {\x}) := \int_{\R^d} W({\x}-{\z}) \varrho_\mathrm{per}({\z}) \, \d {\z},
\]
with $\varrho_\mathrm{per}$ the periodic extension of $\varrho$ to $\R^d$.

Expanding the derivative in the right-hand side of~\eqref{eq:FP} gives
\[
   \part{\varrho}{t} =  \frac12 \, \grad\cdot\Big( \varrho \, \grad G(W\ast\varrho)^2 + G(W\ast\varrho)^2 \, \grad \varrho \Big).
\]
We observe that the convection term $\grad\cdot\left( \varrho \, \grad G(W\ast\varrho)^2 \right)$ 
induces the eventual formation of aggregates, competing against the smoothing action of the purely diffusive term $\grad\cdot \left(G(W\ast\varrho)^2 \, \grad \varrho \right)$.
When suitable model parameters are chosen, stationary patterns emerge as equilibria that balance these two mechanisms.

A linearized stability analysis of a given constant steady state $\varrho \equiv \varrho_0 > 0$
has been carried out in~\cite{BHW:2012}. It leads to the system of ODEs
\(   \label{Fourier}
    \part{\widehat{\varrho}_n}{t} + |q_n|^2 \frac{G(\varrho_0)}{2}
         \left( G(\varrho_0) + 2G'(\varrho_0)\varrho_0 \widehat{W}_n\right) \widehat{\varrho}_n = 0 \,,
\)
for the Fourier transformed perturbation $\tilde\varrho = \tilde\varrho(t, \x)$ of $\varrho_0$,
\(  \label{def:Fourier}
   \widehat \varrho_n : = \int_{\Td} \tilde\varrho(\x) e^{-i q_n\cdot \x}\,\d \x, \qquad
   \widehat W_n : = \int_{\Td} W(\x) e^{-i q_n\cdot \x}\,\d \x,
\)
with $q_n := 2\pi n$, $n\in \Z^d$. We note that due to the radial symmetry of $W$, all Fourier coefficients $\widehat W_n$ are real.
The zero mode $n=0$ in \eqref{Fourier} is neutral, reflecting conservation of  the total mass.
With the assumption $G(\varrho_0) > 0$ and $G'(\varrho_0) < 0$,
those modes $n\neq 0$ are linearly unstable for which
\(   \label{eq:stability}
    \widehat{W}_n > - \frac{G(\varrho_0)}{2G'(\varrho_0)\varrho_0} \geq 0 \,.
\)
Since $W\in L^1(\Td)$, we have $\widehat{W}_n \to 0$ as $|n|\to+\infty$ by the Riemann--Lebesgue lemma.
Therefore, all modes with $|n|$ larger than a certain threshold are linearly asymptotically stable.

\section{Gradient flow structure}\label{sec:GradientFlow}

In this section we show that, for the exponential response function \eqref{eq:Gexp},
the Fokker--Planck equation~\eqref{eq:FP} possesses a natural gradient flow
structure. Throughout this section, all calculations are formal and we assume
that the solution $\varrho=\varrho(t,x)$ is sufficiently smooth and strictly positive.

Let $M>0$ denote the conserved total mass
\begin{equation}
   \label{eq:mass}
   M:= \int_{\Td} \varrho(t,\x)\,\d\x.
\end{equation}
For nonnegative densities $\varrho=\varrho(x)$ with total mass $M>0$, we define the free energy
\begin{equation}
   \label{eq:freeEnergy}
   \mathcal F_\beta(\varrho)
   := \int_{\Td} \varrho(\log\varrho-1)\,\d\x
      -\beta \int_{\Td} \varrho(W\ast\varrho)\,\d\x.
\end{equation}
The first term in~\eqref{eq:freeEnergy} is the Boltzmann entropy, while the
second term represents an attractive nonlocal interaction. 

\begin{proposition}
Let $G$ be given by~\eqref{eq:Gexp}. Then the Fokker--Planck
equation~\eqref{eq:FP} admits the gradient flow form
\begin{equation}
   \label{eq:gradientFlow}
   \part{\varrho}{t}
   =\grad\cdot\left(
      \mm(\varrho)\grad\frac{\delta\mathcal F_\beta}{\delta\varrho}
      \right),
\end{equation}
with the nonlocal mobility
\begin{equation}
   \label{eq:mobility}
   \mm(\varrho)
   :=\frac12G(W\ast\varrho)^2\varrho.
\end{equation}
In particular, $\mathcal F_\beta$ is nonincreasing along every smooth positive
solution of~\eqref{eq:FP}.
\end{proposition}

\begproof
Using the symmetry of $W$, we calculate the first variation
\begin{equation}
   \label{eq:chemicalPotential}
   \frac{\delta\mathcal F_\beta}{\delta\varrho}
   =\log\varrho-2\beta W\ast\varrho, 
\end{equation}
so that~\eqref{eq:Gexp} gives
\[
   \mm(\varrho)\grad\frac{\delta\mathcal F_\beta}{\delta\varrho}
   =\frac12G(W\ast\varrho)^2\grad\varrho
     -\beta G(W\ast\varrho)^2\varrho\grad(W\ast\varrho)
   =\frac12\grad\left(G(W\ast\varrho)^2\varrho\right).
\]
Taking the divergence, we recover~\eqref{eq:FP}. Finally, multiplying
\eqref{eq:gradientFlow} by $\frac{\delta\mathcal F_\beta}{\delta\varrho}$
and integrating by parts over the torus gives
\begin{equation}
   \label{eq:energyDissipation}
   \frac{\d}{\d t}\mathcal F_\beta(\varrho(t))
   = - \int_{\Td} \mm(\varrho) \left| \grad\frac{\delta\mathcal F_\beta}{\delta\varrho} \right|^2 \d x
   \leq 0.
\end{equation}
\endproof

We note that the representation~\eqref{eq:gradientFlow} is understood as a generalized
Wasserstein-type gradient flow, with the symmetric positive Onsager operator
$\mathcal K[\varrho]\xi =-\grad\cdot\bigl(\mm(\varrho)\grad\xi\bigr)$.
Since the mobility~\eqref{eq:mobility} depends nonlocally on $\varrho$, it differs
from the classical $2$-Wasserstein gradient flow, whose mobility is proportional to $\varrho$.

\subsection{Why the exponential response function is necessary}

We next explain why we require the response function to be of the exponential form~\eqref{eq:Gexp}.
Let us therefore assume $G=G(s)$ to be an arbitrary positive differentiable function.
For even more generality, we replace the logarithmic entropy function $s \mapsto s(\log s - 1)$
in the definition \eqref{eq:freeEnergy} of the free energy by some convex, twice differentiable function $f=f(s)$.
I.e., we consider the free energy of the form
\begin{equation}  \label{def:E}
   \mathcal E_\beta(\varrho) = \int_{\Td} f(\varrho) \,\d\x  -\beta \int_{\Td} \varrho(W\ast\varrho)\,\d\x.
\end{equation}
Its first variation is $\frac{\delta\mathcal E_\beta}{\delta\varrho} = f'(\varrho) -2\beta W\ast\varrho$.
Set $\vartheta:=W\ast\varrho$ and suppose that the Fokker-Planck equation \eqref{eq:FP}
admits the generalized gradient-flow representation
\[
     \part{\varrho}{t}  =  \grad\cdot\left( \mm(\varrho,\vartheta) \grad\frac{\delta\mathcal E_\beta}{\delta\varrho}  \right),
\]
with some positive scalar mobility $\mm = \mm(\varrho,\vartheta)$. The gradient flow flux is
$\mm(\varrho,\vartheta) \left(f''(\varrho)\grad\varrho - 2\beta\grad\vartheta \right)$.
On the other hand, the flux in \eqref{eq:FP} is
\[
   \frac12\grad\left(G(\vartheta)^2\varrho\right) = \frac12 G(\vartheta)^2 \grad\varrho + \varrho G(\vartheta)G'(\vartheta)\grad \vartheta.
\]
Assuming that $\grad\varrho$ and $\grad\vartheta$ can be varied independently, equality of the fluxes requires
\[
    \mm(\varrho,\vartheta) f''(\varrho) = \frac12 G(\vartheta)^2 \qquad\mbox{and}\qquad
    -2\beta \mm(\varrho,\vartheta) = \varrho G(\vartheta)G'(\vartheta).
\]
Eliminating the mobility gives
\[
   \varrho f''(\varrho) = -\beta\frac{G(\vartheta)}{G'(\vartheta)}.
\]
The left-hand side depends only on $\varrho$, whereas the right-hand side depends only on $\vartheta$.
Hence both must be equal to some constant $c\in\R$,
\[
   \varrho f''(\varrho)=c \qquad\mbox{and}\qquad -\beta\frac{G(\vartheta)}{G'(\vartheta)}=c.
\]
Consequently,
\[
   f(\varrho)=c\varrho(\log\varrho-1)+a\varrho+b \qquad\mbox{and}\qquad
   G(\vartheta) = k \exp\left(-\frac{\beta}{c} \vartheta\right),
\]
for some $k>0$, while the integration constants $a$, $b\in\R$ can be set zero
since the affine term $a\varrho+b$ is irrelevant in $\mathcal E_\beta(\varrho)$ under the fixed-mass constraint.
Convexity of $f$ requires $c>0$, and therefore $G$ is decreasing when $\beta>0$. The corresponding mobility is
\[
   \mm(\varrho,\vartheta)=\frac{1}{2c}\varrho G(\vartheta)^2 = \frac{k^2}{2c} \varrho \exp\left(-\frac{2\beta}{c} \vartheta\right).
\]
Obviously, the factor $\frac{k^2}{c}>0$ only changes the timescale in \eqref{eq:FP}, and, therefore, can be set to unity
without loss of generality.
We therefore observe that replacing the Boltzmann entropy by a general convex local entropy does not
allow for more general response functions than exponentials. It only allows a rescaling between the coefficient
of the entropy and the exponential decay rate of $G$.

It is insightful to also consider the following alternative argument:
The Fokker-Planck equation \eqref{eq:FP} always admits the factorization
\(   \label{eq:naturalFactorization}
   \part{\varrho}{t}  = \frac12 \grad\cdot\left({Q(\varrho)} \grad\log Q(\varrho)\right)
   \qquad\mbox{with } Q(\varrho):=\varrho G(W\ast\varrho)^2.
\)
The natural candidate for the chemical potential is therefore $\mu(\varrho) = \log\varrho + 2\log G(W\ast\varrho)$.
For $\mu$ to be the first variation of an energy, its Fr\'echet derivative must be symmetric. Denoting
$g(s):=\frac{G'(s)}{G(s)}$, we have
\[
   D\mu(\varrho)[\eta] = \frac{\eta}{\varrho}+2g(W\ast\varrho)W\ast\eta.
\]
Symmetry of the Fr\'echet derivative therefore requires
\[
   g\bigl(W\ast\varrho(\x)\bigr)W(\x-\z) = g\bigl(W\ast\varrho(\z)\bigr)W(\x-\z).
\]
For a generic positive interaction kernel and arbitrary positive densities $\varrho$, this forces $g$ to be constant,
i.e., $G(s) = k\exp(-\beta s)$.

\begin{remark}
The preceding argument builds upon the natural Onsager factorization
\eqref{eq:naturalFactorization}, or, equivalently,
a free energy of the form \eqref{def:E},
composed of entropy and a symmetric pair interaction term.
It does not exclude the possibility that, for a particular nonexponential $G$,
one could introduce a different metric or additional state
variables and obtain a gradient flow formulation of the problem.
However, the author has not found any such structure
(apart from the trivial case $W\equiv\mbox{const}$ when \eqref{eq:FP}
reduces to the heat equation with a constant diffusion coefficient).

Of particular interest would be the power-law case $G(s) = s^{-\alpha}$ with $\alpha>0$.
Here we only note that there is a variational structure for the stationary equation, where
$W\ast\varrho=\lambda\varrho^{p-1}$ for some $\lambda>0$,
where we denoted $p:=1+ \frac{1}{2\alpha}$.
This is the Euler–Lagrange equation for the nonlinear Rayleigh quotient
\[
   \frac{\displaystyle\int_{\Td}\varrho(W*\varrho)\,\d x}
     {\displaystyle \left(\int_{\Td}\varrho^p\,\d x\right)^{2/p}}.
\]
Therefore, although the evolution itself does not appear to be the gradient flow of the Rayleigh quotient,
all positive stationary states admit an exact variational characterization.
\end{remark}

\subsection{Relation to the classical McKean--Vlasov equation}
\label{subsec:McKeanVlasovComparison}
Observe that the free energy~\eqref{eq:freeEnergy} is the standard attractive
McKean--Vlasov free energy \cite{Carrillo-Gvalani-Pavliotis-Schlichting:2020, Carrillo-McCann-Villani:2003}.
The corresponding classical Wasserstein gradient
flow, also known as the granular media equation, is
\begin{equation}
   \label{eq:McKeanVlasovGF}
   \part{\varrho}{t}
   =\grad\cdot\left(\varrho\grad
   \frac{\delta\mathcal F_\beta}{\delta\varrho}\right)
   =\laplace\varrho-2\beta\grad\cdot
   \bigl(\varrho\grad(W\ast\varrho)\bigr).
\end{equation}
By contrast, the spontaneous aggregation model \eqref{eq:FP} takes the form
\begin{equation}
   \label{eq:spontaneousVersusMcKeanVlasov}
   \part{\varrho}{t}
   = \frac{1}{2} \grad\cdot\left( \varrho \exp(-2\beta W\ast\varrho) \grad\frac{\delta\mathcal F_\beta}{\delta\varrho}\right).
\end{equation}
Thus the two equations have the same energy but different Onsager operators:
the mobility in \eqref{eq:McKeanVlasovGF} is $\mm(\varrho)=\varrho$, while
the mobility in \eqref{eq:spontaneousVersusMcKeanVlasov} is $\mm(\varrho)=\frac{1}{2}\varrho e^{-2\beta W\ast\varrho}$.

This observation has several immediate consequences. On a connected domain,
the positive stationary states of both models are precisely the constrained
critical points of $\mathcal F_\beta$. Therefore, the two models have the
same stationary equation, global minimizers, 
and energetic classification of stationary states. Their transient
dynamics, spectra, and convergence rates need not coincide, because these
depend on the mobility. In particular, the mobility in the spontaneous
aggregation model is exponentially suppressed in regions where
$W\ast\varrho$ is large. Rearrangement of mass inside a dense cluster can
therefore occur on a much slower time scale than in the classical
McKean--Vlasov dynamics. This provides a plausible explanation for slow
coarsening and appearance of metastable multi-cluster patterns,
as observed in the numerical simulations in \cite{BHW:2012, EH:2026, HO:2015}.

\subsection{Convexity of the free energy}\label{sec:convex}
We give a sufficient condition for convexity of $\mathcal F_\beta$.
Obviously, the entropy part $\int_{\Td} \varrho (\log\varrho - 1) \,\d x$ is strictly convex in $\varrho$.
To study the convexity of the interaction term
\[
   \mathcal I_\beta(\varrho) := -\beta \int_{\Td} \varrho(W\ast\varrho)\,\d\x,
\]
we calculate its second variation along admissible tangent directions $\eta$,
which are the zero-mean functions $\int_{\Td} \eta\,\d x = 0$. We readily have
\[
   \left. \totk{}{\epsilon}{2} \mathcal I_\beta(\varrho+\epsilon\eta) \right|_{\epsilon=0}
       = -2\beta \int_{\Td} \eta(W\ast\eta)\,\d\x.
\]
Hence convexity of $\mathcal I_\beta$ is equivalent to nonpositivity of the convolution quadratic form on
the subspace of functions with vanishing mean. This can be further characterized using Parseval’s identity,
\[
   \int_{\Td}\eta(W\ast\eta)\,\d x =\sum_{n\in\mathbb Z^d} \widehat W_n|\widehat\eta_n|^2,
\]
with the Fourier coefficients $\widehat W_n$ given by \eqref{def:Fourier}.
Since $\widehat\eta_0=0$, we conclude that $\mathcal I_\beta$ is convex if and only if
$\widehat W_n \leq 0$ for all $n\neq 0$. This is also a sufficient condition for strict convexity
of the energy functional $\mathcal F_\beta$.

\subsection{Existence of global minimizers}\label{subsec:global}

For a fixed total mass $M>0$ we introduce the admissible class
\begin{equation*}
   \mathcal A_M  :=\left\{\varrho\in L^1(\Td):
      \varrho\geq0,\quad  \int_{\Td}\varrho\,\d\x=M,\quad \int_{\Td}\varrho\log\varrho\,\d\x<+\infty  \right\}.
\end{equation*}

\begin{proposition}
\label{prop:existenceGlobalMinimizer}
Let $W\in C(\Td)$ be radially symmetric. Then the problem
\begin{equation}
   \label{eq:variationalProblem}
   \inf_{\varrho\in\mathcal A_M}\mathcal F_\beta[\varrho]
\end{equation}
admits a global minimizer $\varrho$.
\end{proposition}

\begproof
We first show that the free energy is bounded from below.
Indeed, for any $\varrho\in\mathcal A_M$ we have
\[
   \int_{\Td} \varrho(W\ast\varrho)\,\d\x \leq \Norm{W}_{L^\infty(\Td)} M^2.
\]
Moreover, Jensen's inequality for the convex function $s\mapsto s(\log s -1)$ gives
\[
   \int_{\Td} \varrho(\log\varrho-1)\,\d\x \geq M \left( \log\frac{M}{|\Td|}- 1 \right).
\]
Consequently,
\begin{equation}
   \label{eq:energyLowerBound}
   \mathcal F_\beta(\varrho)
   \geq M\log\frac{M}{|\Td|}-M
      -\beta\|W\|_{L^\infty(\Td)}M^2.
\end{equation}
and the infimum in~\eqref{eq:variationalProblem} is finite.

If $(\varrho_n)_{n\in\mathbb N}$ is a minimizing sequence, the entropy
bound and the de la Vall\'{e}e--Poussin criterion give uniform integrability.
Therefore, after an eventual extraction of a subsequence, $\varrho_n$ converges weakly in
$L^1(\Td)$ to some $\varrho\in\mathcal A_M$. The entropy
is weakly lower semicontinuous. Moreover, the continuity of $W$, together
with compactness of $\Td$, implies that
$W\ast\varrho_n$ converges uniformly to $W\ast\varrho$.
Consequently,
\[
   \int_{\Td}\varrho_n(W\ast\varrho_n)\,\d\x  \to   \int_{\Td}\varrho(W\ast\varrho)\,\d\x,
\]
and the direct method of calculus of variations \cite{Dacorogna:2008} gives the existence of a minimizer.

\endproof

\section{Stationary states and their variational characterization}
\label{sec:StationaryStates}

We now use the gradient flow structure derived in Section~\ref{sec:GradientFlow}
to study the stationary solutions of~\eqref{eq:FP}. Throughout this section,
$G$ is given by~\eqref{eq:Gexp}, the total mass is fixed according
to~\eqref{eq:mass}, and we set
\begin{equation}
   \label{eq:homogeneousDensity}
   \varrho_0:=\frac{M}{|\Td|}.
\end{equation}
Obviously, the constant density $\varrho_0$ is a
stationary solution. We first characterize all positive stationary solutions of~\eqref{eq:FP}.

\begin{proposition}
\label{prop:stationaryCharacterization}
Let $\varrho$ be a smooth positive density satisfying~\eqref{eq:mass}. Then the
following statements are equivalent:
\begin{enumerate}[label={\rm(\roman*)}]
   \item $\varrho$ is a stationary solution of~\eqref{eq:FP};
   \item $\varrho$ is a critical point of $\mathcal F_\beta$ under the mass
   constraint~\eqref{eq:mass};
   \item there exists a constant $\lambda\in\mathbb R$ such that
   \begin{equation}  \label{eq:minimizerEL}
   \log\varrho-2\beta W\ast\varrho=\lambda;
   \end{equation}
   \item $\varrho$ satisfies the Gibbs-type equation
   \begin{equation}
   \label{eq:minimizerFixedPoint}
   \varrho(\x)  =M\frac{\exp\bigl(2\beta(W\ast\varrho)(\x)\bigr)}
   {\displaystyle\int_{\Td}
   \exp\bigl(2\beta(W\ast\varrho)(\z)\bigr)\,\d\z}.
\end{equation}
\end{enumerate}
\end{proposition}

\begproof
Let us denote $Q(\varrho):=G(W\ast\varrho)^2\varrho$.
If $\varrho$ is stationary, then $\laplace Q(\varrho)=0$ on the torus.
Consequently, $Q(\varrho)$ is constant. In view
of~\eqref{eq:Gexp}, this is equivalent to
$\varrho\exp(-2\beta W\ast\varrho)=\mathrm{const}$,
which gives~\eqref{eq:minimizerEL}. Conversely,
\eqref{eq:minimizerEL} implies that $\log Q(\varrho)$ is constant and hence that
$\varrho$ is stationary.

By~\eqref{eq:chemicalPotential}, equation~\eqref{eq:minimizerEL} is
the Euler--Lagrange equation for $\mathcal F_\beta$ under the fixed-mass
constraint \eqref{eq:mass}, with $\lambda$ the corresponding Lagrange multiplier. 
Exponentiating~\eqref{eq:minimizerEL} gives
\begin{equation*}
   \varrho(\x)=e^\lambda  \exp\bigl(2\beta(W\ast\varrho)(\x)\bigr).
\end{equation*}
Integration over $\Td$ and use of~\eqref{eq:mass} determine $e^\lambda$ and
give~\eqref{eq:minimizerFixedPoint}.
\endproof

The nonlinear Gibbs-type self-consistency equation~\eqref{eq:minimizerFixedPoint}
gives the following qualitative information about the stationary profiles.

\begin{proposition}
\label{prop:stationaryBounds}
Let $\varrho$ be a positive stationary solution of \eqref{eq:FP} with mass $M>0$.
If $W\in C^k(\Td)$, then $\varrho\in C^k(\Td)$. Moreover,
\begin{equation}
   \label{eq:logGradientBound}
   \|\grad\log\varrho\|_{L^\infty(\Td)}
   \leq 2\beta M\|\grad W\|_{L^\infty(\Td)}
\end{equation}
whenever $W\in W^{1,\infty}(\Td)$, and
\begin{equation}
   \label{eq:contrastBound}
   \frac{\max_{\Td}\varrho}{\min_{\Td}\varrho}
   \leq\exp\bigl(2\beta M\operatorname{osc}_{\Td} W\bigr),
\end{equation}
where
$\operatorname{osc}_{\Td} W:=\max_{\Td} W-\min_{\Td} W$.
The Lagrange multiplier in~\eqref{eq:minimizerEL} satisfies
\begin{equation}
   \label{eq:lambdaBound}
   \lambda\leq\log\varrho_0-2\beta\varrho_0,
\end{equation}
with $\varrho_0$ given by \eqref{eq:homogeneousDensity}.
Equality in \eqref{eq:lambdaBound} holds if and only if $\varrho\equiv\varrho_0$.
\end{proposition}

\begproof
The regularity statement follows directly from~\eqref{eq:minimizerFixedPoint}.
Differentiating~\eqref{eq:minimizerEL} gives
\begin{equation*}
   \grad\log\varrho=2\beta(\grad W)\ast\varrho,
\end{equation*}
which implies~\eqref{eq:logGradientBound}. Moreover, for every $\x,\z\in\Td$,
\begin{equation*}
   \log\frac{\varrho(\x)}{\varrho(\z)}
   = 2\beta\bigl((W\ast\varrho)(\x)-(W\ast\varrho)(\z)\bigr),
\end{equation*}
which implies
\[
   \frac{\max_{x\in\Td} \varrho(\x)}{\min_{z\in\Td} \varrho(\z)} 
   = \exp\left( 2\beta\left( \max_{x\in\Td}(W\ast\varrho)(\x) - \min_{z\in\Td}(W\ast\varrho)(\z)\right) \right)
   = \exp\left( 2\beta\operatorname{osc}_{\Td}(W\ast\varrho) \right).
\]
Since $\varrho$ is nonnegative and has mass $M$, we have
$\operatorname{osc}_{\Td}(W\ast\varrho) \leq M\operatorname{osc}_{\Td} W$,
and~\eqref{eq:contrastBound} follows.

Finally,~\eqref{eq:minimizerEL} combined with \eqref{eq:minimizerFixedPoint} gives
\begin{equation}
   \label{eq:lambdaFormula}
   \lambda=\log M-\log\left(
   \int_{\Td}\exp\bigl(2\beta W\ast\varrho\bigr)\,\d\x\right).
\end{equation}
By~\eqref{eq:Wnorm} and~\eqref{eq:mass},
\begin{equation*}
   \frac1{|\Td|}\int_{\Td} W\ast\varrho\,\d\x =\frac{M}{|\Td|}=\varrho_0.
\end{equation*}
Jensen's inequality therefore gives
\begin{equation*}
   \frac1{|\Td|}\int_{\Td}
   \exp\bigl(2\beta W\ast\varrho\bigr)\,\d\x
   \geq\exp(2\beta\varrho_0).
\end{equation*}
Inserting this into~\eqref{eq:lambdaFormula} proves~\eqref{eq:lambdaBound}.
Equality in Jensen's inequality takes place if and only if $W\ast\varrho$ is
constant. Equation~\eqref{eq:minimizerFixedPoint} then implies that
$\varrho$ is constant and hence equal to $\varrho_0$.
\endproof

We next derive an identity that provides a necessary condition for the
existence of spatially inhomogeneous stationary states.

\begin{proposition}
\label{prop:steadyStateIdentity}
Let $\varrho$ be a positive stationary solution with mass $M>0$ and $\varrho_0$ be given by \eqref{eq:homogeneousDensity}.
Set $r:=\varrho-\varrho_0$.
Then
\begin{equation}
   \label{eq:steadyStateIdentity}
   \int_{\Td}(\varrho-\varrho_0)
   \log\frac{\varrho}{\varrho_0}\,\d\x
   =2\beta\int_{\Td} r(W\ast r)\,\d\x.
\end{equation}
In particular, every nonconstant stationary state satisfies
\begin{equation}
   \label{eq:positiveInteraction}
   \int_{\Td} r(W\ast r)\,\d\x>0.
\end{equation}
\end{proposition}

\begproof
The homogeneous state \eqref{eq:homogeneousDensity} satisfies
\begin{equation*}
   \log\varrho_0-2\beta W\ast\varrho_0
   =\log\varrho_0-2\beta\varrho_0=:\lambda_0.
\end{equation*}
Subtracting this identity from~\eqref{eq:minimizerEL} gives
\begin{equation*}
   \log\frac{\varrho}{\varrho_0}-2\beta W\ast r
   =\lambda-\lambda_0.
\end{equation*}
Multiplication by $r$ and integration over $\Td$,
using $\int_{\Td} r\,\d\x=0$, gives \eqref{eq:steadyStateIdentity}.

The integrand on the left-hand side of \eqref{eq:steadyStateIdentity}
is pointwise nonnegative, since $(a-b)(\log a-\log b)\geq 0$ for all $a$, $b>0$,
and vanishes identically only if $\varrho\equiv\varrho_0$.
Then \eqref{eq:positiveInteraction} immediately follows.
\endproof

Since $W$ is even, Parseval's identity gives
\begin{equation}
   \label{eq:Parseval}
   \int_{\Td} r(W\ast r)\,\d\x
   =\sum_{n\in\mathbb Z^d\setminus\{0\}}
   \widehat W_n|\widehat r_n|^2.
\end{equation}
Consequently, a nonconstant stationary state can exist only if at least one
nonzero Fourier coefficient $\widehat W_n$ is positive.
In other words, if $\widehat W_n\leq0$ for all $n\in\mathbb Z^d\setminus\{0\}$,
then $\varrho\equiv\varrho_0$ is the only stationary state with mass $M$.
Recall that, according to Section \ref{sec:convex}, in this case the energy functional $\mathcal F_\beta$
is strictly convex. Then, obviously, $\varrho_0$ is its unique minimizer on the set of functions with mass $M$.

Moreover, the second variation of the free energy at $\varrho_0$ on the tangent space of the mass constraint is
\begin{equation*}
   \delta^2\mathcal F_\beta(\varrho_0)[\eta,\eta]
   =\int_{\Td}\frac{\eta^2}{\varrho_0}\,\d\x
   -2\beta\int_{\Td}\eta(W\ast\eta)\,\d\x,
   \qquad \int_{\Td}\eta\,\d\x=0.
\end{equation*}
Using \eqref{eq:Parseval}, we have
\begin{align*}
   \delta^2\mathcal F_\beta(\varrho_0)[\eta,\eta]
   &=\sum_{n\in\mathbb Z^d\setminus\{0\}}
   \left(\frac1{\varrho_0}-2\beta\widehat W_n\right)
   |\widehat\eta_n|^2.
\end{align*}
Consequently, $\varrho_0$ loses local minimality in the mode $n$ precisely
when
\begin{equation*}
   \widehat W_n>\frac1{2\beta\varrho_0}.
\end{equation*}
Since $G'(s)=-\beta G(s)$, this is exactly the linear instability
condition~\eqref{eq:stability}. 
Thus the thresholds for linear instability and loss of local
minimality of the homogeneous state coincide.

Finally, we derive a simple sufficient condition under which the homogeneous
density~\eqref{eq:homogeneousDensity} is the unique global minimizer of the free energy.

\begin{proposition}
\label{prop:homogeneousGlobalMinimizer}
Let $W\in L^\infty(\Td)$ be even and let
\[
   \operatorname{ess\,osc}_{\Td} W := 
   \operatorname*{ess\,sup}_{\Td} W -  \operatorname*{ess\,inf}_{\Td} W.
\]
If
\begin{equation}
   \label{eq:smallBetaGlobalMinimizer}
   M \beta \operatorname{ess\,osc}_{\Td} W < 1,
\end{equation}
then $\varrho_0$ is the unique global minimizer of $\mathcal F_\beta$ among
nonnegative densities of mass $M$. 
\end{proposition}

\begproof
Let $\varrho\in\mathcal A_M$ satisfy $\mathcal F_\beta(\varrho)\leq\mathcal F_\beta(\varrho_0)$,
and set $r:=\varrho-\varrho_0$.
Adding a constant to $W$ does not change the quadratic form on zero-mean
functions. Therefore, for every $c\in\mathbb R$,
\(
   \int_{\Td} r(W\ast r)\,\d\x
   =\int_{\Td}\int_{\Td}
   \bigl(W(\x-\z)-c\bigr)r(\x)r(\z)\,\d\z\,\d\x
   \leq \|W-c\|_{L^\infty(\Td)}\|r\|_{L^1(\Td)}^2.
   \label{eq:interactionPinskerStep}
\)
Note that
\[
   \operatorname{ess\,osc}_{\Td} W = 2 \inf_{c\in\mathbb R}\|W-c\|_{L^\infty(\Td)}.
\]
Consequently, taking the infimum over $c\in\R$ in \eqref{eq:interactionPinskerStep} gives
\[
      \int_{\Td} r(W\ast r)\,\d\x \leq \frac12 \|r\|_{L^1(\Td)}^2 \operatorname{ess\,osc}_{\Td} W.
\]
The Csisz\'{a}r--Kullback--Pinsker inequality \cite{Csiszar:1967}, applied after normalizing both
$\varrho$ and $\varrho_0$ by their common mass $M$, yields
\[
   \|r\|_{L^1(\Td)}^2  = \Norm{\varrho-\varrho_0}_{L^1(\Td)}^2
   \leq  2M \int_{\Td}\varrho\log\frac{\varrho}{\varrho_0}\,\d\x.
\]
The relative entropy identity
\[
   \int_{\Td}\varrho\log\frac{\varrho}{\varrho_0}\,\d\x
   =  \mathcal F_\beta(\varrho)-\mathcal F_\beta(\varrho_0) + \beta\int_{\Td} r(W\ast r)\,\d\x
\]
using $\mathcal F_\beta(\varrho) \leq \mathcal F_\beta(\varrho_0)$,
then readily gives
\[
   \|r\|_{L^1(\Td)}^2 \leq  M\beta \|r\|_{L^1(\Td)}^2 \operatorname{ess\,osc}_{\Td} W.
\]
Consequently, if \eqref{eq:smallBetaGlobalMinimizer} holds, then necessarily $\|r\|_{L^1(\Td)}=0$.
\endproof


\section{Single-cluster structure of global minimizers in 1D}\label{sec:SingleCluster}

The variational formulation provides additional information about global
minimizers that is not available for arbitrary stationary states.
In particular, we use symmetric rearrangement to show that a global
minimizer in \eqref{eq:variationalProblem} can be chosen to consist of a
single cluster and, under the additional strict monotonicity assumption
on $W$, that this property holds for every global minimizer;
see, e.g., \cite{Carrillo2019, GGHPS:2025} for related results.
The proof relies on the one-dimensional periodic Riesz rearrangement
inequality \cite{Baernstein,Christ,LiebLoss}.

Let us identify $\mathbb T^1$ with $(-1/2,1/2]$
and recall that the sampling kernel $W(x)$ is supported in $(-1/2,1/2)$,
even and decreasing in the distance $|x|$.
Denote $W^\#$ its symmetric decreasing rearrangement on $\mathbb T^1$, centered at the origin.

\begin{proposition}
\label{prop:singleClusterMinimizer}
Let $d=1$ and let $W\in C(\mathbb T^1)$ satisfy $W=W^\#$.
The variational problem~\eqref{eq:variationalProblem}
admits a global minimizer $\varrho\in\mathcal A_M$ which is, up to a translation, even and nonincreasing,
i.e., for some $x_0\in\mathbb T^1$,
\(   \label{eq:propVarrho}
   \varrho(\x)=\varrho^\#(\x-\x_0).
\)
If $W=W^\#$ is, in addition, strictly decreasing on $[0,1/2]$, then \eqref{eq:propVarrho} holds for
every global minimizer of~\eqref{eq:variationalProblem}.
\end{proposition}

\begproof
Let $\varrho\in\mathcal A_M$ be a global minimizer of~\eqref{eq:variationalProblem},
whose existence is provided by Proposition \ref{prop:existenceGlobalMinimizer}.
Denote by $\varrho^\#$ its symmetric decreasing rearrangement on $\mathbb T^1$, centered at the origin.
Thus $\varrho^\#$ is even and nonincreasing on $[0,1/2]$, and it is
equimeasurable with $\varrho$. In particular,
\begin{equation}
   \label{eq:entropyRearrangement}
   \int_{\mathbb T^1}\Phi(\varrho^\#)\,\d\x
   =\int_{\mathbb T^1}\Phi(\varrho)\,\d\x
\end{equation}
for every Borel function $\Phi$ for which the integrals are well defined.
Taking $\Phi(s)=s\log s$, we see that rearrangement leaves the entropy
term in $\mathcal F_\beta$ unchanged.

Moreover, since $W^\#=W$, the periodic Riesz rearrangement inequality gives
\begin{equation}
   \label{eq:periodicRiesz}
   \int_{\mathbb T^1}\varrho(W\ast\varrho)\,\d\x
   \leq  \int_{\mathbb T^1}\varrho^\# (W\ast\varrho^\#)\,\d\x.
\end{equation}
Since the interaction term appears in~\eqref{eq:freeEnergy} with a negative sign,
\eqref{eq:entropyRearrangement} and~\eqref{eq:periodicRiesz} imply
\begin{equation}
   \label{eq:energyRearrangement}
   \mathcal F_\beta[\varrho^\#] \leq\mathcal F_\beta[\varrho].
\end{equation}
Since $\varrho\in\mathcal A_M$ is a global minimizer of~\eqref{eq:variationalProblem},
then $\varrho^\#\in\mathcal A_M$ is also a global minimizer. This proves the
existence of a symmetric decreasing minimizer.

Moreover, let us observe that since the rearrangement preserves the entropy
term in $\mathcal F_\beta$ due to \eqref{eq:entropyRearrangement},
$\mathcal F_\beta[\varrho^\#] = \mathcal F_\beta[\varrho]$ implies that
\[
   \int_{\mathbb T^1}\varrho(W\ast\varrho)\,\d\x
   = \int_{\mathbb T^1}\varrho^\#  (W\ast\varrho^\#)\,\d\x.
\]
Under the assumption that $W=W^\#$ is strictly decreasing on $[0,1/2]$,
this directly implies \eqref{eq:propVarrho} for some $x_0\in\mathbb T^1$,
see, e.g., \cite{Burchard, Christ}.
Hence the property holds for every global minimizer of~\eqref{eq:variationalProblem}.
\endproof

We emphasize that Proposition~\ref{prop:singleClusterMinimizer} concerns
global minimizers, and not all stationary states.
In particular, it does not exclude that stationary states of \eqref{eq:FP} may
possess multi-cluster critical points.

\section{Long-time behavior}
\label{sec:LongTime}
Another direct consequence of the gradient flow structure of \eqref{eq:FP} is
a characterization of the long time behavior of its solutions.
We only give the formal main idea of the approach, which is based on the observation
that a \mbox{\L}ojasiewicz--Simon inequality is available for \eqref{eq:FP},
due to its similarity to the McKean-Vlasov equation mentioned in Section \ref{subsec:McKeanVlasovComparison}.
Development of a complete rigorous theory is beyond the scope of this short paper
and will be the subject of a future work.

Let $\varrho(t)$ be a smooth positive global solution of~\eqref{eq:FP}. We
write the energy dissipation identity~\eqref{eq:energyDissipation} as
\begin{equation}
   \label{eq:dissipationNotation}
   \frac{\d}{\d t}\mathcal F_\beta(\varrho(t))
   =-\mathcal D_\beta(\varrho(t)),
\end{equation}
where
\begin{equation*}
   \mathcal D_\beta(\varrho)
   :=\frac12\int_{\Td} G(W\ast\varrho)^2\varrho
   \left|\grad \mu_\beta(\varrho)  \right|^2   \,\d\x,
\end{equation*}
with the chemical potential $\mu_\beta(\varrho) := \log\varrho-2\beta W*\varrho$.
The lower bound~\eqref{eq:energyLowerBound} implies that
$\mathcal F_\beta(\varrho(t)) \to \mathcal F_\infty$ as $t\to+\infty$, and, moreover,
\begin{equation}
   \label{eq:finiteTotalDissipation}
   \int_0^\infty\mathcal D_\beta(\varrho(t))\,\d t
   =\mathcal F_\beta(\varrho(0))-\mathcal F_\infty<+\infty.
\end{equation}
Under suitable compactness and continuity assumptions, a LaSalle-type
invariance argument implies that the $\omega$-limit set of the trajectory
is nonempty and contained in the set of stationary states; see
\cite{Hale:1988} and, for an invariance principle in the space of
probability measures, \cite[Theorems~2.12 and~2.13]{Carrillo-Gvalani-Wu:2023}.
Consequently, the distance of $\varrho(t)$ from the set of stationary
states converges to zero in the topology in which the trajectory is
precompact.

A fundamental observation is that the energy dissipation $\mathcal D_\beta$ is uniformly equivalent to
\begin{equation*}
   \widetilde{\mathcal D}_\beta(\varrho)
   :=  \int_{\Td} \varrho  \left|\grad \mu_\beta(\varrho) \right|^2  \,\d\x,
\end{equation*}
which is the dissipation rate of the McKean-Vlasov equation \eqref{eq:McKeanVlasovGF}.
Indeed, since any nonnegative solution is uniformly bounded in $L^1(\Td)$ by the
conservation of mass \eqref{eq:mass}, and since $W\in L^\infty(\Td)$, we have
\[
   \|W*\varrho\|_{L^\infty}\leq B:= M \Norm{W}_{L^\infty(\Td)}.
\]
Then we readily have
\(  \label{eq:DD}
   \frac12e^{-2 \beta B} \widetilde{\mathcal D}_\beta(\varrho)
   \leq
    \mathcal D_\beta(\varrho)
   \leq
    \frac12e^{2 \beta B} \widetilde{\mathcal D}_\beta(\varrho).
\)
This facilitates the application of the \mbox{\L}ojasiewicz--Simon inequality, recently established in \cite{Choi:2026}.
It requires the following quantitative analyticity of $W$,
\[
   \left| \grad^k W(x) \right| \leq \frac{C_W^{k+1} k!}{|x|^{d-2+k}} \qquad\mbox{for all }  x\in\Td\setminus\{0\} \mbox{ and  all } k\in\N,
\]
for some constant $C_W>0$.
Moreover, as it is formulated for probability densities, we introduce the rescaled quantities
\[
   \widetilde\varrho := \varrho/M, \qquad \widetilde\varrho_\infty := \varrho_\infty/M, \qquad \widetilde W := -2\beta M W.
\]
Then, up to the additive constant $M(\log M-1)$,
\[
   \mathcal F_\beta(\varrho) = M \left( \int_{\Td} \widetilde\varrho \log \widetilde\varrho\,\d x
      + \frac12 \int_{\Td} \widetilde\varrho(\widetilde W\ast \widetilde\varrho)\,\d x \right),
\]
and $\grad\bigl(\log \widetilde\varrho+\widetilde W\ast \widetilde\varrho\bigr) = \grad\mu_\beta(\varrho)$.
Then \cite[Theorem~1.1]{Choi:2026} states that if $\widetilde\varrho_\infty$ is a positive stationary state,
then there exist $C$, $\sigma>0$ and $\theta\in [1/2,1)$ such that
\begin{equation*}
   \left|\mathcal F_\beta(\varrho)-\mathcal F_\beta(\varrho_\infty)\right|^{\theta}
   \leq C \Norm{\grad\mu_\beta}_{L^2_\varrho(\Td)} = C \widetilde{\mathcal D}_\beta(\varrho)^{1/2}
\end{equation*}
whenever $\left\|\widetilde\varrho-\widetilde\varrho_\infty \right\|_{W^{1,p}(\Td)}<\sigma$, where $p>d$.
Combined with \eqref{eq:DD}, we have, after an obvious redefinition of the constant $C>0$,
\(
   \label{eq:LSinequality}
   \left|\mathcal F_\beta(\varrho)-\mathcal F_\beta(\varrho_\infty)\right|^{\theta}
   \leq C {\mathcal D}_\beta(\varrho)^{1/2}. 
\)

Let us now choose $\varrho_\infty$ to be an $\omega$-limit point of the trajectory $\varrho=\varrho(t)$.
Then $\varrho_\infty$ is stationary and the continuity of the energy in the respective topology gives
$\mathcal F_\beta(\varrho_\infty)=\mathcal F_\infty$.
In particular,
\[
   e(t):=\mathcal F_\beta(\varrho(t)) -\mathcal F_\beta(\varrho_\infty) \geq 0.
\]
Once the trajectory $\varrho=\varrho(t)$ has entered the neighborhood in which
\eqref{eq:LSinequality} holds, the energy dissipation identity \eqref{eq:dissipationNotation} gives
\begin{equation*}
   \tot{}{t} e(t) = -\mathcal D_\beta(\varrho(t)) \leq-C^{-2}e(t)^{2\theta}.
\end{equation*}
Integration yields algebraic (if $1/2 < \theta < 1$) or exponential (if $\theta=1/2$) decay of $e=e(t)$ to zero.

To obtain convergence of the solution $\varrho=\varrho(t)$, one uses the
\mbox{\L}ojasiewicz--Simon inequality once more. Indeed,
\[
   -\frac{\d}{\d t}e(t)^{1-\theta}
   = (1-\theta)e(t)^{-\theta} \, \mathcal D_\beta(\varrho(t))
   \geq \frac{1-\theta}{C}\mathcal D_\beta(\varrho(t))^{1/2}.
\]
Consequently,
\begin{equation}
   \label{eq:finiteLength}
   \int_{t_0}^{\infty}
   \mathcal D_\beta(\varrho(t))^{1/2}\,\d t<+\infty
\end{equation}
whenever the trajectory remains in the neighborhood in which
\eqref{eq:LSinequality} holds.
We write \eqref{eq:FP} as the continuity equation
\[
   \partial_t\varrho+\div(\varrho v)=0,
   \qquad
   v=-\frac12G(W\ast\varrho)^2\grad\mu_\beta(\varrho).
\]
By the characterization of absolutely continuous curves in Wasserstein space through the continuity equation,
see \cite[Theorem~8.3.1]{AGS:2008}, we have for the 2-Wasserstein metric derivative 
\[
   | \dot{\widetilde\varrho} |_{W_2}(t)
   \leq \left(\int_{\Td} |v(t,x)|^2 \widetilde\varrho(t,x)\,\d x\right)^{1/2}
   \leq \frac{e^{\beta B}}{\sqrt{2M}}\,
   \mathcal D_\beta(\varrho(t))^{1/2}.
\]
Consequently, for $t>s\geq t_0$,
\[
   W_2(\widetilde\varrho(s),\widetilde\varrho(t))
   \leq \frac{e^{\beta B}}{\sqrt{2M}}
   \int_s^t\mathcal D_\beta(\varrho(r))^{1/2}\,\d r.
\]
Thus \eqref{eq:finiteLength} implies that
$\widetilde\varrho(t)$ is Cauchy in the complete space of Borel probability measures on $\Td$
equipped with the $2$-Wasserstein distance.
It therefore converges to a probability measure.
Since $\widetilde\varrho_\infty$ is an $\omega$-limit point of
the trajectory, this limit must coincide with $\widetilde\varrho_\infty$.

A rigorous implementation of the argument 
would require a globally positive solution,
uniform-in-time bounds ensuring parabolic regularity and precompactness of
the trajectory in $W^{1,p}(\Td)$ with $p>d$.
Under these assumptions, one first extracts a
stationary $\omega$-limit point $\varrho_\infty$. A standard trapping
argument then shows that, once the solution enters a sufficiently small
neighborhood of $\varrho_\infty$, it cannot leave it, and the above
finite-length estimate yields $\varrho(t) \to \varrho_\infty$ as $t\to+\infty$
in the Wasserstein distance.
Additional uniform regularity estimates may then be used to upgrade
the result to convergence in a strong Sobolev or classical
topology.
We refer to the proof of Theorem~1.3 in \cite{Choi:2026} for details.


\end{document}